\documentstyle{amsppt}
\magnification=\magstep1
\topmatter
\title
The Distribution of the Eigenvalues of Hecke Operators
\\
\endtitle
\author
J.B.~Conrey \\
W.~Duke \\
D.W.~Farmer
\endauthor
\thanks Research of the first and second author supported in 
part by a grant from the NSF.
\endthanks
\address
 Department of Mathematics, Oklahoma State University,
Stillwater, OK 74078
\endaddress
\address
 Department of Mathematics, Rutgers University,
New Brunswick, NJ 08903
\endaddress
\abstract
For each prime $p$, we determine
the distribution of the $p^{{th}}$ Fourier coefficients of the 
Hecke eigenforms of large weight for the full modular group. 
As $p\to\infty$, this distribution tends to the Sato--Tate distribution.
\endabstract
\endtopmatter
\document

\def \frac #1#2{{ #1\over #2 }}
\def\intl{\int\limits}

\def\({\left(}
\def\){\right)}

\magnification=\magstep1

\NoBlackBoxes

\head
1.  Introduction and statement of results
\endhead

Let
$$
f(z)=\sum_{n=1}^\infty a_{f,n} e(nz)
$$
be a cusp form of weight $k$ for the full modular group. We assume
that $f(z)$ is an eigenfunction of all the Hecke operators~$T_n$,
and we set $a_f(n)=n^{1-k\over 2} a_{f,n}$.
The assumption that $f(z)$ is a Hecke eigenform gives:
$$
\eqalignno{ a_f(mn)&=a_f(n)a_f(m) &\hbox{if } (m,n)=1\ \ \ \ \ \cr
&\cr
a_f(p^j)&=a_f(p)a_f(p^{j-1}) - a_f(p^{j-2}) &\hbox{for $p$ prime}\ \ \ \ \ \ \ \cr
&\cr
|a_f(n)|&\le d(n),
}$$
where $d(n)$ is the divisor function.
The above relations were first conjectured by Ramanujan based on a 
computation of the coefficients of the weight 12 cusp form
$\Delta(z)=\sum \tau(n)e(nz)$.  The two equations were proven for
$\tau(n)$ by Mordell, using what are now known as the Hecke operators.
The inequality was proven by Deligne as a consequence of his
proof of the Weil conjectures.  Those results determine everything
about $a_f(n)$ except for the distribution of the $a_f(p)\in[-2,2]$.

Define $\theta_f(p)\in[0,\pi]$ by $a_f(p)=2\cos\theta_f(p)$.  It is  conjectured
that for each $f$ the $\theta_f(p)$ are uniformly distributed with respect 
to the
Sato--Tate measure 
$$
{2\over \pi} \sin^2\theta\, d\theta.
$$
This conjecture was first made for 
cusp forms associated to
non--CM elliptic curves, and was 
extended to Hecke eigenforms for the full modular group by 
Serre~\cite{Ser}.
This conjecture appears to be quite deep, for it has been 
shown~\cite{Ogg}~\cite{Mur} to be equivalent to the nonvanishing on
the line $Re(s)=1$ of all the $m^{{th}}$ symmetric power
$L$--functions $L_m(s)$ associated to $f(z)$.
 
In this paper we change perspective slightly by fixing $p$ and
looking at the distribution of $a_f(p)$ as $f(z)$ varies.  That is,
we are looking at the distribution of the eigenvalues of $T_p$.

We write $S_k(\Gamma(1))$ for the space of cuspforms of weight~$k$
for the full modular group, and we write $f\in S_k$ to mean
$f\in S_k(\Gamma(1))$ and $f$ is a Hecke eigenform.
Following Selberg, we write $\sigma_k(T_n)$ for 
the trace of the Hecke operator
$T_n$ acting on $S_k(\Gamma(1))$.  Because of our normalization,
$\sigma_k(T_n)=n^{(k-1)/2} \sum_{f\in S_k}a_f(n)$.

\proclaim{Theorem 1}  Suppose $p$ is prime.  As $k\to\infty$, 
the set $\{\theta_f(p):\ f\in S_k\}$ 
becomes uniformly distributed with respect to the measure
$$
{2\over \pi}\(1+{1\over p}\) {\sin^2\theta \over
	{\(1-{1\over p}\)^2} + {4\over p}\sin^2\theta}\,d\theta .
$$
\endproclaim

The same distribution was found by Sarnak~\cite{Sar} for the
$p^{th}$ coefficient of Maass forms averaged over the Laplacian
eigenvalues.
Sarnak has pointed out to us that the above measure is
the $p$--adic Plancherel measure, and it is also the spectral
measure of the nearest--neighbor Laplacian on a
$p+1$ regular tree (see~\cite{LPS}).

\proclaim{Theorem 2}  As $p\to\infty$, with $p$ prime,  and $k\to\infty$ with
$k>e^p$, the set $\{\theta_f(p):\ f\in S_k\}$ becomes
uniformly distributed with respect to the Sato--Tate measure.
\endproclaim

A result similar to this was obtained by Birch~\cite{B}, who
considered the distribution of $\theta(p)$ over the elliptic
curves over {\bf F}$_p$.  Similar results have also been found
for Kloosterman sums~\cite{K}~\cite{A}~\cite{Mic} and cubic exponential
sums~\cite{L}.  In each case, results have been obtained 
by fixing the prime $p$ and varying the other parameters in the problem.
The original problem of varying the prime, keeping the
other parameters fixed, remains inaccessible. 

The idea behind the proofs of the Theorems is as follows.  
Using the Selberg trace formula
we obtain a simple expression for $\sum_{f\in S_k} a_f(p^j)$, 
valid for large $k$. 
Using the recurrence
relation satisfied by $a_f(p^j)$ leads to an expression for
$\sum_{f\in S_k} a_f(p)^j$.  Such an expression for all $j$ is sufficient to
determine the distribution of $\{a_f(p):\ f\in S_k\}$.

\head
2.  Initial Lemmas
\endhead

Our starting point is the Selberg trace formula applied to 
the Hecke operators.  The following is formula (4.5) in~\cite{Sel}.

\proclaim{Lemma 1}  We have
$$
\eqalign{\sigma_k(T_n)=&-{1\over2}\sum_{-2\sqrt{n}<m<2\sqrt{n}}
	H(4n-m^2){\eta_m^{k-1} - {\overline{\eta}}^{k-1}_m \over
		\eta_m- {\overline{\eta}}_m} \cr
	&-{\sum_{{d|n}\atop d\le \sqrt{n}}}' d^{k-1} 
		+ \delta(\sqrt{n}){k-1\over12}n^{{k\over2} -1},
}$$ 
where $H(d)$ is the class number, $\delta(x)=1$ if $x$ is an integer
and is 0 otherwise, $\sum^\prime$ means that a term with $d=\sqrt{n}$
is counted with weight 1/2, and
$$
\eta_m = {m+i(4n-m^2)^{1\over2}\over 2}.
$$
\endproclaim

The above formula directly leads to an expression for 
$\sum_{f\in S_k} a_f(p^j)$.
Using the estimate $H(d)\ll d^{{1\over 2}+\varepsilon}$ gives

\proclaim{Lemma 2}  Suppose $p$ is prime.  As $k\to\infty$,
$$\eqalignno{
\sum_{f\in S_k} a_f(p^j)&={k\over 12} {1\over p^{j\over2}} 
+ O(p^{{j\over 2}+\varepsilon})
&\hbox{$j$ even\ \ \ \ \ \ \ \ \ \ \ \ \ }\cr
&\cr
&= O(p^{{j\over 2}+\varepsilon}) &j \hbox{ odd\ \ \ \ \ \ \ \ \ \ \ \ \ }
}$$
\endproclaim

To make use of the above Lemma, it is necessary to relate 
$a_f(p)^n$ to $a_f(p^j)$.  This involves the 
Chebyshev polynomials $U_n(x)$.  We have

\proclaim{Lemma 3}  Suppose $p$ is prime and write $a_f(p)=2\cos \theta_f(p) $.  Then
$$
a_f(p^j)=U_j(\cos \theta_f(p))
$$
and
$$
a_f(p)^n = \sum_{j=0}^n h_n(j) a_f(p^j),
$$
where
$$
h_n(j)={2^{n+1}\over \pi}\intl_0^\pi \cos^n\theta\, \sin  (j+1)\theta  \,
			\sin \theta \,d\theta .
$$
\endproclaim

\demo {Proof}  The first expression follows directly from the
recurrence relation $a_f(p^j)=a_f(p)a_f(p^{j-1}) - a_f(p^{j-2})$ and the
fact that the $U_j(x)$ satisfy essentially the same recurrence.  
To prove the second, note that the $h_m(j)$ satisfy
$$
x^n = \sum_{j=0}^n h_n(j) U_n(x/2).
$$
Since the $U_n(x)$ are orthonormal on $[-1,1]$ with
respect to the measure ${2\over\pi}\sqrt{1-x^2}\,dx$, we have
$$
h_n(j)={2^{n+1}\over \pi}\intl_{-1}^1 x^n U_n(x) \sqrt{1-x^2}\,dx .
$$
Changing variables $x\mapsto \cos\theta$ gives the stated formula.
\enddemo

\head
3.  Proof of the Theorems
\endhead

Both Theorems follow immediately from this Proposition.

\proclaim{Proposition 1}  Suppose $p$ is prime.  As $k\to\infty$, 
$$
\sum_{f\in S_k} a_f(p)^n = 
{k\over 6\pi}\(1+{1\over p}\) \intl_0^\pi 2^n \cos^n \theta\,
{\sin^2\theta \over
	1+{1\over p^2}-{2\over p}\cos 2\theta}\,d\theta 
+ O\(p^{{n\over 2}+\varepsilon}\) 
$$
\endproclaim

\demo{Proof}  We assume $n$ is even, the case of
odd $n$ being identical.  Starting with the second expression
in Lemma 3, and then applying Lemma 2 gives,
$$
\eqalign{
\sum_{f\in S_k} a_f(p)^n&=\sum_{f\in S_k} \sum_{j=0}^n h_n(j)a_f(p^j)\cr
&=
{k\over 12}   \sum_{j=0}^{n/2} {h_n(2j)\over p^j} 
	+ O\(p^{{n\over 2}+\varepsilon}\)\cr
&={k\over 12}   \sum_{j=0}^\infty {h_n(2j)\over p^j} 
+ O\(p^{{n\over 2}+\varepsilon}\),
}$$
the completion of the sum on $j$ being valid because $h_n(j)=0$
for $j>n$.  Now inserting the expression for $h_n(j)$ from Lemma~3,
switching the sum and integral, and then evaluating the sum on $j$ gives,
$$
\eqalign{
\sum_{f\in S_k} a_f(p)^n&={k\over 6\pi}\sum_{j=0}^\infty {1\over p^j}
\intl_0^\pi 2^n\cos^n\theta\, \sin  (2j+1)\theta  \,
                        \sin \theta \,d\theta 
	+ O\(p^{{n\over 2}+\varepsilon}\)\cr
&={k\over 6\pi}\intl_0^\pi 2^n\cos^n\theta\, \sin \theta 
\sum_{j=0}^\infty {\sin  (2j+1)\theta \over p^j} \,d\theta 
+ O\(p^{{n\over 2}+\varepsilon}\)\cr
&=
{k\over 6\pi}\(1+{1\over p}\) \intl_0^\pi 2^n \cos^n \theta\,
{\sin^2\theta \over
        1+{1\over p^2}-{2\over p}\cos 2\theta}\,d\theta
+ O\(p^{{n\over 2}+\varepsilon}\), 
}$$
as claimed.
\enddemo

The deduction of the Theorems from Proposition~1 is by Weyl's
criterion: the sequence $(b_n)$ is uniformly distributed on $[a,\,b]$
with respect to the measure $d\mu$ provided
$$
\lim_{N\to\infty} {1\over N}\sum_{n=1}^N g(b_n)
=
\intl_a^b g\,d\mu,
$$
for all $g\in C([a,\,b])$.  Proposition~1 gives the necessary equality
for $g(x)=x^n$, uniformly in $k$, which is sufficient because the 
polynomials are dense in $C([a,\,b])$. This completes the proofs of
the Theorems.

\Refs

\item{[A]} A. Adolphson, {\it On the distribution of angles of 
Kloosterman sums}, J. Reine Angew. Math. {\bf 395}, (1989).

\item{[B]} B. J. Birch, {\it How the number of points of an elliptic 
curve over a fixed prime field varies},  J. London Math. Soc. {\bf 43} (1968).

\item{[K]} N. M. Katz, {\it Gauss sums, Kloosterman sums, 
and monodromy groups}, Ann. of Math. Stud.,
116, Princeton, 1988. 

\item{[L]} R. Livn\'e, {\it The average distribution of cubic exponential sums},
J. Reine Angew. Math. {\bf 375/376} (1987).

\item{[LPS]} A. Lubotzky, R. Phillips, and P Sarnak, 
{\it Ramanujan graphs}, Combinatorica {\bf 8} no.3 (1988)

\item{[Mic]} P. Michel, {\it Autour de la conjecture de Sato-Tate pour
les sommes de Kloosterman. I}, Invent. Math. {\bf 121} (1995).

\item{[Mur]} V. K. Murty, {\it On the Sato--Tate conjecture},
in {\it Number theory related to Fermat's last theorem},
Cambridge, MA, 1981.

\item{[Ogg]}  A.~P.~Ogg, {\it A remark on the Sato--Tate conjecture},
Invent.~Math. {\bf 9}, (1970).

\item{[Sar]} P.~Sarnak, {\it  Statistical properties of eigenvalues 
of the Hecke operators}, in {\it Analytic number theory
and Diophantine problems (Stillwater, OK, 1984)}, Progr. Math., 70,
Birkhauser, Boston, 1987. 

\item{[Sel]} A.~Selberg, {\it Harmonic analysis and
discontinuous groups}, in {\it  Collected papers. Vol. I},
Springer--Verlag, Berlin, 1989.

\item{[Ser]} J.P. Serre,  {\it Abelian l--adic representations
and elliptic curves},  Benjamin, New York, 1968.

\endRefs

\enddocument